\documentclass[11pt,twoside]{amsart}
\usepackage{amssymb}
\usepackage{mathrsfs}
\def\flex#1{\mathrel{\mathop{\kern 0pt\hbox to 
10mm{\rightarrowfill}}\limits_{#1
\rightarrow \infty}}}

\def\Flex#1{\mathrel{\mathop{\kern 0pt\hbox to 
10mm{\rightarrowfill}}\limits_{#1
\rightarrow \infty}}}

\font\bb=msbm10 at 12pt

\def\P{\hbox{\bb P}}

\def\n {\bf n}
\def\1{1\hspace{-1.2mm}\mbox{{\normalsize I}}}

\newtheorem{Th}{Theorem}

\newtheorem{Lemme}[Th]{Lemma}

\newtheorem{Prop}[Th]{Proposition}
\newcommand{\R}{\mathbb{R}}

\newcommand{\F}{\mathcal{F}}

\title{Brownian motion in the quadrant with oblique repulsion from the sides}
\begin{document}
\maketitle
\centerline{Dominique L\'{e}pingle\footnote{Universit\'{e} d'Orl\'{e}ans, MAPMO-FDP, F-45067 Orl\'{e}ans 
(dominique.lepingle@univ-orleans.fr)}}

\begin{abstract}
We consider the problem of strong existence and uniqueness of a Brownian motion  forced to
stay in the nonnegative quadrant by an electrostatic repulsion from the sides that works obliquely. When 
the direction of repulsion is normal, the question has previously been solved with the help of convex analysis.  
To construct the solution, we start from the normal
 case and then we use as  main tool a
comparison lemma. The results are reminiscent of the study of a  Brownian motion with oblique reflection 
in a wedge. Actually, the same skew symmetry condition is involved when looking for a stationary distribution in product form. 
The terms of the product are now gamma distributions in place of exponential ones.
An associated purely deterministic problem is also considered.
\end{abstract} 

\section{Introduction}

In the late seventies the study of heavy traffic limits in open multi-station queueing networks has put the question of existence 
and properties of the Brownian motion obliquely reflected on the sides of a wedge and more generally on the faces of a polyhedron. 
In the following decade there was an extensive literary output on that topic, among which we mention the works of 
Harrison, Reiman, Williams and collaborators (\cite{H78},\cite{HR81},\cite{VW85},\cite{W87},\cite{HW87},\cite{DW96}, 
to cite a few of them). But there is another way than normal or oblique reflection to prevent Brownian motion from overstepping 
 a linear barrier. We may add as drift term the gradient of a concave function that explodes in the neighborhood 
of the faces of the polyhedron. To be more specific, let $\n_1,\ldots,\n_k$ be unit vectors in $\R^d$ and $b_1,\ldots,b_k$ be 
real numbers. The state space $S$ is defined by 
\[
  S:= \{{\bf x}\in \R^d: {\bf n}_r.{\bf x}\geq b_r, r=1,\ldots,k\}\,.
\]
Let $\phi_1,\ldots,\phi_k$ be $k$ convex $C^1$ functions on $(0,\infty)$ with $\phi_r(0+)=+\infty$ for any $r=1,\ldots,k$. 
The potential function $\Phi$ on $S$ is defined by
\[
  \Phi({\bf x}):= \sum_{r=1}^k\phi_r({\bf n}_r.{\bf x}-b_r)\,.
\]
From the general existence and uniqueness theorem on multivalued stochastic differential systems established in 
\cite{C95} and \cite{C98}, completed with identification of the drift term as in Lemma 3.4 in
\cite{CL97},
we know there exists a unique strong solution living in $S$ 
to the equation
\begin{equation}
 \label{eq:convexe}
 \begin{array}{rll}
  d{\bf X}_t & = & d{\bf B}_t - \nabla\Phi({\bf X}_t)dt \\
   {\bf X}_0 & \in & S
 \end{array}
\end{equation}
where ${\bf B}$ is a Brownian motion in $\R^d$. As was proved in Proposition 4.1 of \cite{L10}, the hypotheses 
$\phi_r(0+)=+\infty$ for any $r=1,\ldots,k$ entail that there is no additional  boundary process of local time type 
in the r.h.s. of (\ref{eq:convexe}) since 
the repulsion forces are sufficiently strong. As
\[
  -\nabla \Phi({\bf x})= -\sum_{r=1}^k\phi_r^{\prime}({\bf n}_r.{\bf x}-b_r){\bf n}_r\,,
\]
the repulsion from the faces $F_r=\{{\bf x}\in S: {\bf n}_r.{\bf x}=b_r\}$ points in a normal direction into the interior of $S$. 
We now introduce  vectors ${\bf q}_1,\ldots,{\bf q}_k$  with ${\bf q}_r.{\bf n}_r=0$ for $r=1,\ldots,k$ and consider a new 
drift function 
\begin{equation}
 \label{eq:general}
  -\sum_{r=1}^k\phi_r^{\prime}({\bf n}_r.{\bf x}-b_r)({\bf n}_r+{\bf q}_r) \,.
\end{equation}
This is a singular drift and convex analysis cannot be used as in \cite{C95} or \cite{C98} to get strong existence and
uniqueness of the associated stochastic differential system. 

In this paper we will concentrate on the nonnegative quadrant in $\R^2$ as state space and a drift term that derives 
from  electrostatic repulsive forces. 
We obtain strong existence and uniqueness for a large set of parameters. Our results are reminiscent of the thorough study in 
\cite{VW85} of the oblique reflection in a wedge. However a key tool in this work was an appropriate harmonic function that made 
a weak approach possible and fruitful. Thus full results were obtained, while our strong approach merely provides a partial 
answer to the crucial question of hitting the corner.

For theoretical as well as practical reasons, a great deal of interest was taken in the question of existence and computation
of the invariant measure of the Brownian motion with a constant drift vector and oblique reflection 
(\cite{HR81},\cite{W87},\cite{HW87}). Under the assumption that the directions of reflection satisfy a skew symmetry condition,
it was proved that the invariant measure has exponential product form density. Motivated by the so-called Atlas model of equity markets
presented in \cite{F02}, 
some authors (\cite{PP08},\cite{IK10},\cite{IP11}) have recently studied Brownian motions on the line with rank dependent 
local characteristics. This model is strongly related to Brownian motion reflected in polyhedral domains. 
The invariant probability density has an explicit exponential product form when the volatility is constant \cite{PP08} and a sum of 
products of exponentials form when the volatility coefficients depend on the rank (\cite{IP11},\cite{IK12}). 
This last kind of density was previously 
obtained in \cite{DM09} for a Brownian motion in a wedge with oblique reflection.

A neighboring way has been recently explored in \cite{OO12}. Here the process is a Brownian motion with a drift term that is continuous 
and depends obliquely,  via a regular potential function, on the position of the process  relative to a polyhedral domain. Under the same
 skew symmetry condition as in \cite{HW87}, the invariant density has an explicit product form again. In  Section 5,
we consider a Brownian motion with a constant drift living in the nonnegative
quadrant with oblique electrostatic repulsion from the sides. 
Under the  skew symmetry condition, there is still an  invariant measure in product form. Now the terms of the product are two gamma
distributions with explicit parameters.

In the last section we consider the same model of oblique electrostatic repulsion, this time without Brownian term. 
We are interested in the question of existence and uniqueness of this deterministic differential system when the starting point is the corner. There is a trivial solution and the proofs of uniqueness of the previous stochasting setting are still in force in this simpler case. Taking advantage of the explicit form of this solution, we can obtain weaker conditions for uniqueness. 

\section{The setting}

The general state space is the nonnegative quadrant $S=\R_+\times \R_+$. 
The corner ${\bf 0}=(0,0)$ will play a crucial role and in some cases
it will be necessary to restrict the state space to the punctured nonnegative quadrant 
$S^{{\bf 0}}=S\setminus \{{\bf 0}\}$.

Let $(B_t,C_t)$ be a Brownian motion in the plane starting from ${\bf 0}$, adapted to a filtration ${\F}=({\F}_t)$ 
with usual conditions. Let $\alpha,\beta,\gamma,\delta$ be four real constants with $\alpha>0,\delta>0$. We say that an
${\F}$-adapted continuous process $(X,Y)$ with values in $S$ is a Brownian motion with electrostatic oblique repulsion
from the sides if for any $t\geq 0$
\begin{equation}
 \label{eq:obli}
 \begin{array}{lllll}
 X_t & = & X_0+B_t+\alpha\int_0^t\frac{ds}{X_s} + \beta\int_0^t\frac{ds}{Y_s} & \geq & 0 \\
 Y_t & = & Y_0+C_t+\gamma\int_0^t\frac{ds}{X_s} + \delta\int_0^t\frac{ds}{Y_s} & \geq & 0
 \end{array}
\end{equation}
where $X_0$ and $Y_0$ are non-negative ${\F}_0$-measurable random variables. Each coordinate $X_t$ or $Y_t$ may
vanish, so to make sense we must have a.s. for any $t\geq 0$
\[
\begin{array}{lll}
  \int_0^t{\bf 1}_{\{X_s=0\}}ds =0 & \quad & \int_0^t{\bf 1}_{\{Y_s=0\}}ds =0\\
  \int_0^t{\bf 1}_{\{X_s>0\}}\frac{ds}{X_s}< \infty  & \quad & \int_0^t{\bf 1}_{\{Y_s>0\}}\frac{ds}{Y_s}< \infty \,.
\end{array}
\]
The drift in (\ref{eq:obli}) is of type (\ref{eq:general}) with $d=k=2$ and
\[
 \begin{array}{ll}
\phi_1(x)=-\alpha\log(x) & \phi_2(y)= -\delta\log(y)  \\
 b_1=0 & b_2=0 \\
 {\bf n}_1=(1,0)  &  {\bf n}_2=(0,1) \\
 {\bf q}_1=(0,\frac{\gamma}{\alpha})  &  {\bf q}_2=(\frac{\beta}{\delta},0) 
 \end{array}
\]
The case with $\beta=\gamma=0$ is a particular case of (\ref{eq:convexe}). In the sequel, we will note $(U,V)$
the solution of the system
\begin{equation}
 \label{eq:normal}
   \begin{array}{lllll}
 U_t & = & X_0+B_t+\alpha\int_0^t\frac{ds}{U_s}  & \geq & 0 \\
 V_t & = & Y_0+C_t + \delta\int_0^t\frac{ds}{V_s} & \geq & 0 \,.
 \end{array}
\end{equation}
The processes $U$ and $V$ are independent Bessel processes (if $X_0$ and $Y_0$ are independent variables).
Actually, $U$ is a Bessel process with index $\alpha-\frac{1}{2}$, and the point $0$ is intanstaneously reflecting for $U$ if 
$\alpha<\frac{1}{2}$ and polar if $\alpha\geq\frac{1}{2}$. Moreover, $U^2+V^2$ is the square of a Bessel process 
with index $\alpha+\delta$, and so the corner ${\bf 0}$ is polar for $(U,V)$ in any case.

Comparison between $X$ and $U$, $Y$ and $V$ will play a key role in the construction of the solution $(X,Y)$ and the 
study of its behavior near the sides of the quadrant. The following simple lemma will be of constant use.
\begin{Lemme}
For $T>0$, $\alpha>0$, let $x_1$ and $x_2$ be nonnegative continuous solutions on $[0,T]$ of the equations
\[
  \label{Th:comp}  
\begin{array}{lll}
   x_1(t) & = & v_1(t) + \alpha\int_0^t\frac{ds}{x_1(s)} \\
   x_2(t) & = & v_2(t) + \alpha\int_0^t\frac{ds}{x_2(s)}
  \end{array}
\]
where  $v_1$, $v_2$ are continuous functions such that $0\leq v_1(0)\leq v_2(0)$, and $v_2-v_1$ is 
nondecreasing.
Then $x_1(t)\leq x_2(t)$ on $[0,T]$.
\end{Lemme}
{\bf Proof}.
Assume there exists $t\in (0,T]$ such that $x_2(t)<x_1(t)$. Set
  \[
   \tau:=\max \{s\leq t : x_1(s)\leq x_2(s)  \} \,.
 \]
Then,
\[
  \begin{array}{lll}
  x_2(t)-x_1(t) &= & x_2(\tau)-x_1(\tau)+(v_2(t)-v_1(t))-(v_2(\tau)-v_1(\tau)) +
  \alpha  \int_{\tau}^t(\frac{1}{x_2(s)}-\frac{1}{x_1(s)})ds \\
   &\geq & 0 \,,
  \end{array}
\]
a contradiction. $\hfill  \blacksquare$

We will also  need the following consequence of the results in \cite{C95} or \cite{C98} on multivalued
stochastic differential systems, completed with the  method used in \cite{L10} to check the lack of additional boundary process.
\begin{Prop}
\label{Th:inter}
Let $\alpha>0$, $\delta\geq 0$, ${\bf \sigma}=(\sigma_j^i; i,j=1,2) $ a $2\times2$-matrix, $(B_1,B_2)$ 
a Brownian motion in the plane, $b_1$ and $b_2$ two Lipschitz functions on $\R^2$,
$Z^1_0$ and $Z_0^2$ two ${\F}_0$-measurable  nonnegative random variables. There exists a unique  solution 
$(Z^1,Z^2)$ to the system
\begin{equation}
 \label{eq:mixte}
 \begin{array}{lll}
  Z^1_t & = & Z^1_0+ \sigma_1^1B^1_t+\sigma_2^1B^2_t + \alpha \int_0^t\frac{ds}{Z^1_s} + \int_0^tb_1(Z^1_s,Z^2_s)ds \\
  Z^2_t & = & Z^2_0+ \sigma_1^2B^1_t+\sigma_2^2B^2_t + \delta\int_0^t\frac{ds}{Z^2_s} + \int_0^tb_2(Z^1_s,Z^2_s)ds 
 \end{array}
\end{equation}
with the conditions $Z^1_t\geq 0$ if $\delta=0$ and $Z^1_t\geq 0,Z^2_t\geq 0$ if $\delta>0$.
\end{Prop}

It is worth noticing that the solutions to (\ref{eq:obli}) enjoy the Brownian scaling property. It means that if $(X,Y)$ is a solution to (\ref{eq:obli})
starting from $(X_0,Y_0)$ with driving Brownian motion $(B_t,C_t)$, then for any $c>0$ the process $(X_t^{\prime}:=c^{-1}X_{c^2t}, Y_t^{\prime}:=
c^{-1}Y_{c^2t}; t\geq 0)$ is a solution to (\ref{eq:obli}) starting from $(c^{-1}X_0,c^{-1}Y_0)$ with driving Brownian motion $(c^{-1}B_{c^2t},
c^{-1}C_{c^2t})$.

\section{Avoiding the corner}

We shall see in the next section that  existence and uniqueness of the solution to (\ref{eq:obli}) are easily obtained as soon 
as the solution process keeps away from the corner. Thus the question of attaining the corner in finite time is of great interest.

\begin{Th}
\label{Th:corn}
Let $(X,Y)$ be a solution to (\ref{eq:obli}) in the interval $[0,\tau]\cap [0,\infty)$ where $\tau$ is a ${\F}$-stopping time. 
 We set 
\[
  \tau^{{\bf 0}}:= \inf \{t\in(0,\tau]\cap (0,\infty):(X_t,Y_t)={\bf 0}\}
\]
with the usual convention $\inf \emptyset =  \infty$.
Then $\P(\tau^{{\bf 0}}<\infty)=0$ if one of the following conditions is satisfied:
\begin{enumerate}
 \item
  $C_1: \beta \geq 0$ and $\gamma \geq 0 $
 \item
  $C_{2a} : \alpha\geq \frac{1}{2}$ and  $\beta \geq 0$
 \item
  $C_{2b} : \delta\geq \frac{1}{2}$ and  $\gamma \geq 0$
 \item
 $C_3 :$ There exist $\lambda>0$ and $\mu>0$ such that
  \begin{itemize}
 \item
   $\lambda \alpha + \mu \gamma \geq 0$  
 \item
  $\lambda\beta+\mu \delta \geq 0$
   \item
   $\lambda(\lambda\alpha+\mu \gamma)+\mu(\lambda\beta +\mu \delta) -\frac{1}{2}(\lambda^2 + \mu^2)
   \geq -2\sqrt{\lambda\mu(\lambda\beta+\mu\delta)(\lambda \alpha+\mu\gamma)}$.
   \end{itemize}
 \end{enumerate}
\end{Th}
{\bf Proof}. Assume $\tau=\infty$ a.s. for simplicity of notation.\\
{\it Condition} $C_1$. From Lemma 1 we get $X_t\geq U_t$, $Y_t\geq V_t$, where $(U,V)$ is the solution to (\ref{eq:normal}),
and we know that ${\bf 0}$ is polar for $(U,V)$.\\
{\it Condition} $C_{2a}$ (resp. $C_{2b}$). From Lemma 1 we get $X_t\geq U_t$ (resp. $Y_t\geq V_t$) and in this case $0$ is polar for 
$U$ (resp. V), so $U_t>0$ (resp. $V_t>0$) for $t>0$. \\
{\it Condition} $C_3$.
For $\epsilon>0$ let 
\[
  \begin{array}{lll}
    \sigma^{\epsilon} & = & {\bf 1}_{\{(X_0,Y_0)={\bf 0}\}}\inf \{t>0: X_t+Y_t\geq \epsilon\} \\
    \tau^{{\bf 0},\epsilon} & = & \inf\{t>\sigma^{\epsilon}: (X_t,Y_t)={\bf 0}\} \,.
\end{array}
\]
As $\epsilon\downarrow 0$, $\sigma^{\epsilon}\downarrow 0$ 
and $ \tau^{0,\epsilon}\downarrow \tau^{{\bf 0}}$.
We set $S_t=\lambda X_t+ \mu Y_t$ for $t\geq 0$, $\lambda>0$ and $\mu>0$. From Ito formula we get 
for $t\in[\sigma^{\epsilon},\tau^{{\bf 0},\epsilon})$
\[
 \begin{array}{ll}
  &\log S_t \\
 = & \log S_{\sigma^{\epsilon}} + \int_{\sigma^{\epsilon}}^t\frac{\lambda dB_s+\mu dC_s}{S_s} + 
  (\lambda\alpha+\mu \gamma)\int_{\sigma^{\epsilon}}^t\frac{ds}{X_sS_s}
  +(\lambda \beta+\mu \delta)\int_{\sigma^{\epsilon}}^t\frac{ds}{Y_sS_s} -\frac{1}{2}
     (\lambda^2 +\mu^2)\int_{\sigma^{\epsilon}}^t\frac{ds}{S_s^2} \\
   = & \log S_{\sigma^{\epsilon}} + M_t + \int_{\sigma^{\epsilon}}^t\frac{P(X_s,Y_s)}{X_sY_sS_s^2}ds
 \end{array}
\]
where $M$ is a continuous local martingale and $P(x,y)$ is the second degree homogeneous polynomial
\[
  P(x,y)= \lambda(\lambda\beta+\mu\delta)x^2+\mu(\lambda\alpha+\mu\gamma)y^2+
  (\lambda(\lambda\alpha+\mu\gamma)+\mu(\lambda\beta+\mu\delta)-\frac{1}{2}(\lambda^2+\mu^2))xy \,.
\]
Condition  $C_3$ is exactly the condition for $P$ being nonnegative in $S$. 
Therefore 
\[
  0\leq \int_{\sigma^{\epsilon}}^t \frac{P(X_s,Y_s)}{X_sY_sS_s^2}ds < \infty
\]
and so 
\[
0\leq \int_{\sigma^{\epsilon}}^{\tau^{{\bf 0},\epsilon}} \frac{P(X_s,Y_s)}{X_sY_sS_s^2}ds \leq \infty \,.
\]
As $t\rightarrow \tau^{{\bf 0},\epsilon}$, the continuous local martingale $M$ either converges to a finite limit or oscillates 
between $+\infty$ and $-\infty$.  It cannot converge to $-\infty$ and thus $S_{\tau^{{\bf 0},\epsilon}}>0$
on $\{\tau^{{\bf 0},\epsilon}<\infty\}$, proving that $\P(\tau^{{\bf 0},\epsilon}<\infty)=0$ and finally
$\P(\tau^{{\bf 0}}<\infty)=0$.$\hfill  \blacksquare$ 

\vspace{0.5cm}
\noindent{\bf Example}. When $\alpha=\delta$ and $|\beta|=|\gamma|$, condition $C_3$ is satisfied (with $\lambda=\mu$) if
\begin{equation}
  \begin{array}{llll}
   \bullet & \beta^2\leq \alpha -\frac{1}{4} & \mbox{ when} &  \beta=-\gamma \\
   \bullet & -\beta \leq \alpha-\frac{1}{4} & \mbox{ when} & \beta=\gamma < 0\,.
  \end{array}
\end{equation}

\vspace{0.5cm}
We may also be interested in  hitting  a single side. Then we set
\begin{equation}
 \begin{array}{lll}
  \tau_X^0 & := & \inf \{t>0: X_t=0\} \\
  \tau_Y^0 & := & \inf \{t>0: Y_t=0\} \,.
 \end{array}
\end{equation}
We already know that $\P(\tau_X^0<\infty)=0$ if $\alpha\geq\frac{1}{2}$ and $\beta\geq 0$.
Conversely we can prove
that $\P(\tau_X^0<\infty)=1$ if $\alpha<\frac{1}{2}$ and $\beta\leq 0$. 
If we know that the 
corner is not hit and $\alpha\geq\frac{1}{2}$, we can get rid of the nonnegativity assumption on $\beta$.

\begin{Prop}
\label{Th:side}
Assume $\P(\tau^{{\bf 0}}<\infty)=0$. If $\alpha\geq \frac{1}{2}$, then $\P(\tau^0_X<\infty)=0$.
\end{Prop}
{\bf Proof}. For $\eta>0$ let
\[
  \begin{array}{lll}
    \theta_X^{\eta} & = & {\bf 1}_{\{X_0=0\}}\inf \{t>0: X_t\geq \eta\} \\
    \tau_X^{0,\eta} & = & \inf\{t>\theta_X^{\eta}: X_t=0\} \,.
\end{array}
\]
As $\eta\downarrow 0$, $\theta^{\eta}_X\downarrow 0$  and $ \tau^{0,\eta}_X\downarrow \tau^0_X$.
For $t\in[\theta^{\eta}_X,\tau^{ 0,\eta}_X)$,
\begin{equation}
\label{eq:bordX}
  \begin{array}{lll}
  \log X_t & = & \log X_{\theta_X^{\eta}} + \int_{\theta_X^{\eta}}^t \frac{dB_s}{X_s} +
    (\alpha-\frac{1}{2})\int_{\theta_X^{\eta}}^t \frac{ds}{X_s^2} + \beta \int_{\theta_X^{\eta}}^t \frac{ds}{X_sY_s}\,.
  \end{array}
\end{equation}
Since $\P(\tau^{{\bf 0}}<\infty)=0$, we have $Y_{\tau_X^{0,\eta}}>0$ on $\{\tau_X^{0,\eta}<\infty\}$. On this set,
\[
  \int_{\theta_X^{\eta}}^{\tau_X^{0,\eta}} \frac{ds}{X_s} <\infty \;, 
 \qquad   \int_{\theta_X^{\eta}}^{\tau_X^{0,\eta}} \frac{ds}{Y_s} <\infty
\]
and $X_s>0$ on $[\theta_X^{\eta}, \tau_X^{0,\eta})$, which proves that
\[
   \beta \int_{\theta_X^{\eta}}^{ \tau_X^{0,\eta}}\frac{ds}{X_sY_s}  > - \infty \,.
\]
As $t\rightarrow \tau_X^{0,\eta}$, the local martingale in the r.h.s. of (\ref{eq:bordX}) cannot converge to $-\infty$. 
This entails that $\P(\tau_X^{0,\eta}<\infty)=0$ and therefore $\P(\tau^0_X<\infty)=0$. $\hfill \blacksquare$

\vspace{0.5cm}

We may use again the  method  in Theorem \ref{Th:corn} and Proposition \ref{Th:side} to learn more about hitting the sides
of the quadrant.

\begin{Prop}
\label{Th:sides}
Assume $\alpha\geq \frac{1}{2}$ and $\delta\geq \frac{1}{2}$. Then $\P(\tau^0_X<\infty)=\P(\tau^0_Y<\infty)=0$ if one of the
following conditions is satisfied:
 \begin{enumerate}
  \item
   $\beta>0$
  \item
   $\gamma>0$
  \item 
   $\beta\gamma \leq (\alpha-\frac{1}{2})(\delta-\frac{1}{2})$.
 \end{enumerate}
\end{Prop}

{\bf Proof}. For $\epsilon>0$ let
 \[
   \begin{array}{lll}
    \rho^{\epsilon} & = & {\bf 1}_{\{X_0Y_0=0\}} \inf\{t>0 : X_tY_t\geq \epsilon\} \\
    \rho^{0,\epsilon} & = & \inf\{t>\rho^{\epsilon}:X_tY_t=0\}\;.
  \end{array}
\]
For $\lambda>0$ and $\mu>0$ we set
\[
  R_t=\lambda \log X_t +\mu \log Y_t.
\]
From Ito formula we get for $t\in [\rho^{\epsilon},\rho^{0,\epsilon})$
\[
  \begin{array}{lll}
   R_t & = & R_{\rho^{\epsilon}} + \int_{\rho^{\epsilon}}^t(\frac{\lambda}{X_s}dB_s+\frac{\mu}{Y_s}dC_s) + 
   \int_{\rho^{\epsilon}}^t [\frac{\lambda (\alpha-\frac{1}{2})}
   {X_s^2} + \frac{\mu (\delta-\frac{1}{2})}{Y_s^2} + \frac{(\lambda\beta+\mu\gamma)}
   {X_sY_s} ]ds\\
    & = & R_{\rho^{\epsilon}} + N_t +\int_{\rho^{\epsilon}}^t \frac{Q(X_s,Y_s)}{X_s^2Y_s^2}ds
  \end{array}
\]
where $N$ is a continuous local martingale and $Q(x,y)$ is the second degreee homogeneous polynomial 
\[
  Q(x,y) = \mu(\delta-\frac{1}{2})x^2 + \lambda(\alpha-\frac{1}{2})y^2 +(\lambda \beta +\mu \gamma) xy.
\]
If we obtain that $Q$ is nonnegative on $S$, then the proof will terminate as in Theorem \ref{Th:corn} under condition $C_3$. 
If $\beta>0$ or $\gamma>0$, we easily find $\lambda>0$ and $\mu>0$ such that $\lambda \beta+ \mu \gamma \geq 0$, and then $Q$ is 
nonnegative on $\R^2$. If now $\beta\leq 0$, $\gamma\leq 0$ and $\beta\gamma\leq (\alpha-\frac{1}{2})(\delta-\frac{1}{2})$, we may take 
$\lambda=-\gamma$, $\mu=- \beta$ and we check that $Q(x,y)$ remains nonnegative on $\R^2$ as well. $\hfill  \blacksquare$

\section{Existence and uniqueness}

We now proceed to  the question of  existence and uniqueness of a global solution to (\ref{eq:obli}). We consider separately 
the three cases: $\beta\geq 0$ and $\gamma \geq 0$, then $\beta> 0$ and $\gamma<0$, then $\beta\leq 0$ and $\gamma<0$.
 
\subsection{Case $\beta\geq 0$ and $\gamma \geq 0$}
This is exactly condition $C_1$.

\begin{Th}
 \label{Th:premier}
Assume $\beta\geq 0$ and $\gamma \geq 0$. 
\begin{enumerate}
 \item
  There is a unique  solution to (\ref{eq:obli}) in $S^{{\bf 0}}$.
 \item
  There is a  solution to (\ref{eq:obli}) in $S$ starting from ${\bf 0}$.
 \item
  If $\alpha \delta\geq \beta \gamma$, there is a unique  solution  to (\ref{eq:obli}) in $S$.
\end{enumerate}
\end{Th}
{\bf Proof}. $1.$ Let $a>0$, $\epsilon>0$ and define for $(x,z)\in \R_+\times \R$
\[
  \psi_{\epsilon}(x,z):= \frac{1}{\max(\gamma x+z,\alpha \epsilon)}\,.
\]
This is a Lipschitz function. From Proposition \ref{Th:inter} we know that the system
\begin{equation}
  \label{eq:constru}
   \begin{array}{lll}
   X_t^{\epsilon} & = & X_0+B_t+\alpha\int_0^t\frac{ds}{X_s^{\epsilon}} +\alpha\beta\int_0^t \psi_{\epsilon}
  (X_s^{\epsilon},Z_s^{\epsilon})ds \geq 0 \\
  Z_t^{\epsilon}& = & -\gamma X_0+\alpha(Y_0+{\bf 1}_{\{Y_0=0\}}a)-\gamma B_t+\alpha C_t +\alpha(\alpha\delta-\beta\gamma)
    \int_0^t \psi_{\epsilon}(X_s^{\epsilon},Z_s^{\epsilon})ds
  \end{array}
\end{equation}
has a unique  solution. Let
\[
  \tau_Y^{\epsilon} := \inf \{t>0:\gamma X_t^{\epsilon}+ Z_t^{\epsilon}<\alpha \epsilon \} \,.
\]
If $0<\eta<\epsilon<a$ we deduce from the uniqueness that $(X^{\epsilon},Z^{\epsilon})$ and
$(X^{\eta},Z^{\eta})$ are identical on $ [0,\tau_Y^{\epsilon}]$. Patching together we can set
\[
 \begin{array}{lll}
  X_t & := & \lim_{\epsilon\rightarrow 0} X_t^{\epsilon} \\
  Y_t & := & \lim_{\epsilon\rightarrow 0} \frac{1}{\alpha}(\gamma X_t^{\epsilon} + Z_t^{\epsilon})
 \end{array}
\]
on $\{Y_0>0\} \times [0,\tau_Y^0)$, where 
\[
  \tau_Y^0 := \lim_{\epsilon \rightarrow 0} \tau_Y^{\epsilon} \,.
\]
On this set, $(X,Y)$ is the unique solution to (\ref{eq:obli}). As we noted in the proof of Theorem 3 with condition $C_1$, we have 
$X_t\geq U_t$ and $Y_t\geq V_t$. Therefore, on $\{Y_0>0\} \cap\{\tau_Y^0<\infty\}$,
\[
  \begin{array}{lll}
    \int_0^{\tau_Y^0} \frac{ds}{X_s} \leq \int_0^{\tau_Y^0} \frac{ds}{U_s} <\infty  &  \mbox{ and }  &
    \int_0^{\tau_Y^0} \frac{ds}{Y_s} \leq \int_0^{\tau_Y^0} \frac{ds}{V_s} <\infty 
  \end{array}
\]
and we can define
\begin{equation}
 \begin{array}{lllll}
  X_{\tau_Y^0} & := & \lim_{t\rightarrow \tau_Y^0} X_t & = & X_0+B_{\tau_Y^0} + \alpha \int_0^{\tau_Y^0}\frac{ds}{X_s}
       +\beta \int_0^{\tau_Y^0}\frac{ds}{Y_s} \\
  Y_{\tau_Y^0} & := & \lim_{t\rightarrow \tau_Y^0} Y_t & = & Y_0+C_{\tau_Y^0} + \gamma \int_0^{\tau_Y^0}\frac{ds}{X_s}
       +\delta \int_0^{\tau_Y^0}\frac{ds}{Y_s} \,.
  \end{array}
\end{equation}
We have  $Y_{\tau_Y^0}=0$ and  as ${\bf 0}$ is polar for $(U,V)$, then $X_{\tau_Y^0}>0$. 
In exactly the same way we can construct a solution 
on $\{Y_0>0\}$ in the interval $[T_1,T_2]$, where $T_1=\tau_Y^0$, $T_2= \inf \{t>T_1:X_t=0\}$. Iterating, we get a solution on 
$\{Y_0>0\} \times [0,\lim_{n\rightarrow\infty} T_n)$
where
\[
  \begin{array}{lll}
   T_{2p} & := & \inf\{t>T_{2p-1} : X_t=0\} \\
    T_{2p+1} & := & \inf\{t>T_{2p} : Y_t=0\}\,.
  \end{array}
\]
On $\{ Y_0>0\} \cap \{\lim_{n\rightarrow \infty}T_n<\infty\}$  
we  set $X_{\lim_{n\rightarrow \infty}T_n} := \lim_{p\rightarrow \infty} X_{T_{2p}}=0$ and 
$Y_{\lim_{n\rightarrow \infty}T_n} :=\lim_{p\rightarrow \infty} Y_{T_{2p+1}}=0$. The polarity of ${\bf 0}$ entails 
this is not possible in finite time and thus 
$\lim_{n\rightarrow \infty} T_n = \infty$. So we have obtained a unique global   solution on $\{Y_0>0\}$. In the same
way we obtain a unique global  solution on $\{X_0>0\}$ and as $\P((X_0,Y_0)= {\bf 0})=0$  the proof is complete.

2. Assume now $X_0=Y_0=0$. Let $(y_n)_{n\geq 1}$ be a sequence of real numbers (strictly) decreasing to $0$. From the above paragraph it follows there exists  
for any $n\geq 1$ a  unique solution $(X^n,Y^n)$ with values in $S^{\bf 0}$ to the system
\[
  \begin{array}{lll}
    X_t^n  & = & B_t+\alpha\int_0^t\frac{ds}{X_s^n} + \beta\int_0^t\frac{ds}{Y_s^n}  \\
    Y_t ^n& = & y_n+C_t+\gamma\int_0^t\frac{ds}{X^n_s} + \delta\int_0^t\frac{ds}{Y^n_s} \,.
 \end{array}
\]
 Let 
\[
  \tau:= \inf \{t>0: X^{n+1}_t<X_t^n\} \,.
\]
 Using  Lemma \ref{Th:comp}  we obtain $Y_t^{n+1} \leq Y_t^n$ on $[0,\tau]$. We note that
$(X_{\tau}^n,Y_{\tau}^n)\in S^{{\bf 0}}$ on $\{\tau<\infty\}$. On 
$\{Y_{\tau}^{n+1}=Y_{\tau}^n\}\cap \{\tau<\infty\}$, since $X_{\tau}^{n+1}=X_{\tau}^n$ and 
 the solution starting at time $\tau$ is unique, it follows 
that $X_t^{n+1}=X_t^n$ and $Y_t^{n+1}=Y_t^n$ on $[\tau, \infty)$.
On $\{Y_{\tau}^{n+1}<Y_{\tau}^n\}\cap \{\tau<\infty\}$, the continuity of solutions at time $\tau$ 
entails there exists $\rho>0$ such that $Y_t^{n+1}\leq Y_t^n$ on $[\tau,\tau+\rho]$. A second application of Lemma 1 proves that $X_t^{n+1}\geq X_t^n$ on 
$[\tau,\tau+\rho]$, a contradiction
to the definition of $\tau$. Therefore $\P(\tau=\infty)=1$. It follows that 
$X_t^{n+1}\geq X_t^n$ and $Y_t^{n+1}\leq Y_t^n$ for any $t\in [0,\infty)$,
and we may define
\[
X_t:=\lim_{n\rightarrow\infty} \uparrow X_t^n  \qquad Y_t:=\lim_{n\rightarrow\infty} \downarrow Y_t^n \,.
\]
As $Y_t^n\geq V_t$ where $(U,V)$ is the solution to (\ref{eq:normal}) with $X_0=Y_0=0$, we have
\[
  \begin{array}{lll}
  X_t& = & B_t+ \alpha \lim_{n\rightarrow \infty} \int_0^t \frac{ds}{X_s^n}+ 
    \beta  \lim_{n\rightarrow \infty} \int_0^t \frac{ds}{Y_s^n} \\
  & =&  B_t + \alpha \int_0^t\frac{ds}{X_s}+ \beta \int_0^t\frac{ds}{Y_s} \\
   &< & \infty
  \end{array}
\]
and also 
\[
 \begin{array}{lll}
  Y_t& = &\lim_{n\rightarrow \infty}y_n +  C_t+ \gamma \lim_{n\rightarrow \infty} \int_0^t \frac{ds}{X_s^n}+ 
    \delta  \lim_{n\rightarrow \infty} \int_0^t \frac{ds}{Y_s^n} \\
  & =&  C_t + \gamma \int_0^t\frac{ds}{X_s}+ \delta \int_0^t\frac{ds}{Y_s} \\
   &< & \infty \,.
  \end{array}
\]

3. Assume finally $\alpha\delta-\beta \gamma \geq 0$. As the conclusion holds true if $\beta=\gamma=0$, we may also 
assume $\beta>0$. Let $(X,Y)$ be the solution to (\ref{eq:obli}) with $X_0=Y_0=0$ 
obtained in the previous paragraph and let $(X^{\prime},Y^{\prime})$ be another solution. Considering $(X^n,Y^n)$ again 
and replacing 
$(X^{n+1},Y^{n+1})$ with $(X^{\prime},Y^{\prime})$, the previous proof works and we finally obtain
$X^{\prime}_t\geq X_t$ and $Y^{\prime}_t\leq Y_t$. Then,
\begin{equation}
 \label{eq:uni} 
 \begin{array}{ll}
 & (\delta(X_t-X^{\prime}_t)-\beta(Y_t-Y_t^{\prime}))^2\\  = &
   2 \int_0^t(\delta(X_s-X^{\prime}_s)-\beta(Y_s-Y_s^{\prime}))(\alpha\delta-\beta\gamma)(\frac{1}{X_s}-\frac{1}{X_s^{\prime}})ds \\
   \leq & 0
  \end{array}
\end{equation}
and thus $X^{\prime}_t=X_t$ and $Y_t^{\prime}=Y_t$, proving  uniqueness. Replacing 
$(\delta,\beta)$ with $(\gamma,\alpha)$  in  equation 
(\ref{eq:uni}) we obtain the same conclusion 
if $\gamma>0$. $\hfill \blacksquare$

\subsection{Case $\beta> 0$ and $\gamma<0$}

\begin{Th}
\label{Th:deux}
Assume $\beta> 0$, $\gamma<0$ and one of the conditions $C_{2a}$ or $C_3$ is satisfied. Then, there
exists a unique  solution to (\ref{eq:obli}) in $S^{{\bf 0}}$.
\end{Th}
{\bf Proof}. 
The proof is similar to the proof of 1 in Theorem \ref{Th:premier}. The only change is that now 
$Y_t\leq V_t$. Therefore, on $\{Y_0>0\} \cap\{\tau_Y^0<\infty\}$,
\[
  \delta \int_0^{\tau_Y^0}\frac{ds}{Y_s}  \leq V_{\tau_Y^0} -Y_0 -C_{\tau_Y^0}-\gamma \int_0^{\tau_Y^0} \frac{ds}{U_s} < \infty
\]
and we can define $X_{\tau_Y^0}$ and  $Y_{\tau_Y^0}$ as previously done. $\hfill \blacksquare$

\subsection{Case $\beta\leq 0$ and $\gamma<0$}

In this case we can give a full answer to the question of existence and uniqueness. Our condition of existence is 
exactly the condition found in \cite{W85} for the reflected Brownian in a wedge being a semimartingale, i.e. there is a convex combination 
of the directions of reflection that points into the wedge from the corner.

\begin{Th}
\label{Th:trois}
Assume $\beta\leq 0$ and $\gamma<0$.
\begin{enumerate}
  \item
   If $\alpha\delta>\beta\gamma$, there exists a unique  solution to (\ref{eq:obli}) in $S$.
  \item
  If $\alpha \delta \leq \beta \gamma$, there does not exist any solution.
\end{enumerate}
\end{Th}
{\bf Proof}.
1. Assume first $\alpha \delta > \beta \gamma$.\\
a) {\it Existence}.  Let $(h_n,n\geq 1)$ be a (strictly) increasing sequence of bounded positive 
nonincreasing Lipschitz functions converging to 
$1/x$ on $(0,\infty)$ and to $+\infty$ on $(-\infty,0]$. For instance we can take 
\[
  \begin{array}{lllll} 
      h_n(x) & = & (1-\frac{1}{n}) \frac{1}{x}  & \quad\mbox{ on } & [\frac{1}{n},\infty) \\
                 & = & n-1   & \quad \mbox{ on } & (-\infty,\frac{1}{n}] \,.
  \end{array}
\]
We consider for each $n\geq 1$ the system
\begin{equation}
 \label{eq:appro}
  \begin{array}{lll}
 X_t ^n& = & X_0+B_t+\alpha\int_0^t\frac{ds}{X_s^n} + \beta\int_0^th_n(Y_s^n)ds  \\
 Y_t^n & = & Y_0+C_t+\gamma\int_0^t h_n(X_s^n)ds + \delta\int_0^t\frac{ds}{Y_s^n} \,.
 \end{array}
\end{equation}
From Proposition \ref{Th:inter} it follows there exists a unique  solution to this system. 
 We set
\[
   \tau:= \inf\{s>0: X_s^{n+1}>X_s^n\}\,.
\]
We have $h_{n+1}(X_t^{n+1})\geq h_n(X_t^n)$ on $[0,\tau]$. A first application of Lemma \ref{Th:comp} shows that 
  $Y_t^{n+1}\leq Y_t^n$ on
$[0,\tau]$. Since
$h_{n+1}(Y_{\tau}^{n+1})> h_n(Y_{\tau}^n)$ on $\{\tau<\infty\}$, 
we deduce from the continuity of solutions that
there exists 
$\rho>0$ such that $h_{n+1}(Y_t^{n+1})\geq h_n(Y_t^n)$ on $[\tau, \tau+\rho]$. 
A second application of  Lemma \ref{Th:comp} shows that  $X_t^{n+1}\leq X_t^n$ on $[\tau, \tau+\rho]$,
a contradiction to the definition of 
$\tau$. Thus $\P(\tau=\infty)=1$ proving that on the whole $[0,\infty)$ we have $X_t^{n+1}\leq X_t^n$ and $Y_t^{n+1}\leq Y_t^n$.
Then we can set for any $t\in [0,\infty)$
\[
    X_t:= \lim_{n\rightarrow \infty}X_t^n \quad \mbox{ and } \quad   Y_t:= \lim_{n\rightarrow \infty}Y_t^n \,.
\]
If $\alpha\delta>\beta \gamma$, there is a convex combination of the directions of repulsion that points into 
the positive
quadrant, i.e. there exist $\lambda>0$ and $\mu>0$ such that $\lambda\alpha+\mu \gamma>0$ and $\mu\delta+\lambda
\beta>0$. For $n\geq 1$ and $t\geq 0$, 
\begin{equation}
 \label{eq:majo}
 \begin{array}{lll}
  \lambda U_t+\mu V_t & \geq & \lambda X_t^n+ \mu Y_t^n \\
    & \geq & \lambda X_0 + \mu Y_0 + \lambda B_t + \mu C_t + (\lambda\alpha+\mu\gamma) \int_0 ^t \frac{ds}{X_s^n}                    
                     + (\mu\delta+\lambda \beta) \int_0 ^t \frac{ds}{Y_s^n} \,.
  \end{array}
\end{equation}
Letting $n \rightarrow \infty$ in (\ref{eq:majo}) we obtain
\[
 \int_0^t \frac{ds}{X_s} <\infty \quad \mbox{and} \quad \int_0^t \frac{ds}{Y_s} <\infty \,.
\]
Then we may let $n$ go to $\infty$ in  (\ref{eq:appro}) proving that $(X,Y)$ is a solution to (\ref{eq:obli}).\\
b) {\it Uniqueness}. Let $(X^{\prime},Y^{\prime})$ be another solution to (\ref{eq:obli}). Replacing $(X^{n+1},Y^{n+1})$ with 
$(X^{\prime},Y^{\prime})$ we follow the above proof to obtain for $t\in[0,\infty)$ and $n\geq 1$ 
\[
   X_t^{\prime} \leq X^n_t \quad  \mbox{and} \quad  Y_t^{\prime} \leq Y^n_t
\]
Letting $n\rightarrow \infty$ we conclude
\[
  X_t^{\prime} \leq X_t \quad  \mbox{and} \quad  Y_t^{\prime} \leq Y_t \,.
\]
With the same $\lambda>0$ and $\mu>0$ as above,
\[
 \begin{array}{ll}
  & (\lambda (X_t-X_t^{\prime}) + \mu (Y_t-Y_t^{\prime})) ^2 \\
  = & 2 \int_0^t (\lambda(X_s-X_s^{\prime})+ \mu(Y_s-Y_s^{\prime}))\,[(\lambda \alpha + \mu \gamma)
 (\frac{1}{X_s}-\frac{1}{X_s^{\prime}})+(\mu \delta+\lambda \beta)(\frac{1}{Y_s}-\frac{1}{Y_s^{\prime}})]ds \\
  \leq & 0
 \end{array}
\]
and therefore $X_t^{\prime}=X_t$, $Y_t^{\prime}=Y_t$. \\
2. If $\alpha \delta\leq \beta \gamma$ there exist $\lambda> 0$ and $\mu> 0$ such that $\lambda\alpha+\mu \gamma\leq 0$
and $\mu\delta+\lambda\beta\leq 0$. Thus, if $(X,Y)$ is a solution to (\ref{eq:obli}),
\[
 0 \leq  \lambda X_t + \mu Y_t \leq \lambda X_0 + \mu Y_0  + \lambda B_t + \mu C_t \,.
\]
This is not possible since the paths of the  Brownian motion 
$(\lambda^2+\mu^2)^{-1/2} (\lambda B_t+\mu C_t)$ are not bounded below. 
So there is no global solution.
$\hfill \blacksquare $

In the following pictures, we display the $x$-repulsion direction vector $r_x=(\alpha,\gamma)$ and the $y$-repulsion direction 
vector $r_y=(\beta,\delta)$ in three illustrative instances.

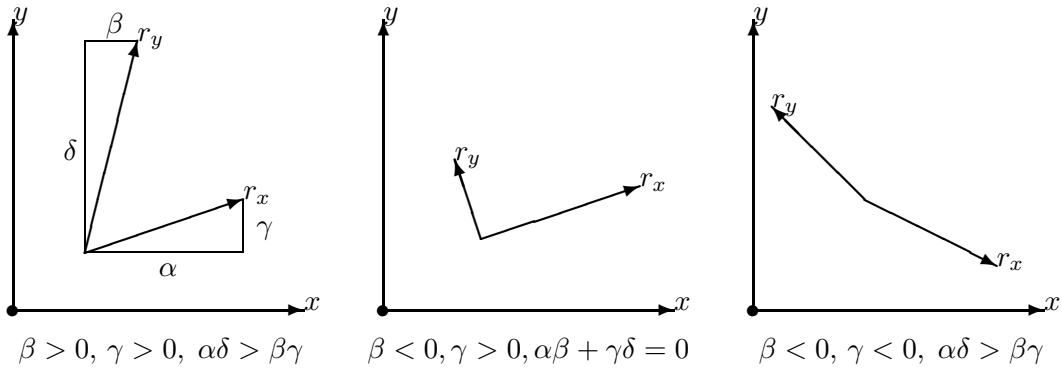
\begin{picture}(400,165)(5,20)
\begin{thicklines}
\put(5,40){$\bullet$}
\put(8,43){\vector(1,0){110}}
\put(8,43){\vector(0,1){110}}
\put(35,65){\vector(3,1){60}}
\put(35,65){\vector(1,4){20}}
\put(95,85){$r_x$}
\put(55,145){$r_y$}
\put(118,43){$x$}
\put(8,153){$y$}
\put(10,25){$\beta>0,\,\gamma>0,\;\alpha\delta> \beta\gamma$}
\put(145,40){$\bullet$}
\put(148,43){\vector(1,0){110}}
\put(148,43){\vector(0,1){110}}
\put(258,43){$x$}
\put(148,153){$y$}
\put(185,70){\vector(3,1){60}}
\put(185,70){\vector(-1,3){10}}
\put(245,90){$r_x$}
\put(175,100){$r_y$}
\put(142,25){$\beta<0,\gamma> 0,\alpha\beta+\gamma\delta=0$}
\put(285,40){$\bullet$}
\put(288,43){\vector(1,0){110}}
\put(288,43){\vector(0,1){110}}
\put(398,43){$x$}
\put(288,153){$y$}
\put(330,85){\vector(2,-1){50}}
\put(330,85){\vector(-1,1){35}}
\put(380,60){$r_x$}
\put(295,120){$r_y$}
\put(290,25){$\beta<0,\,\gamma<0,\;\alpha\delta>\beta\gamma$}

\end{thicklines}
\begin{thinlines}
\put(35,65){\line(0,1){80}}
\put(35,145){\line(1,0){20}}
\put(35,65){\line(1,0){60}}
\put(95,65){\line(0,1){20}}
\put(63,56){$\alpha$}
\put(100,72){$\gamma$}
\put(27,100){$\delta$}
\put(43,148){$\beta$}

\end{thinlines}
\end{picture}

\vspace{1cm}

\section{Product form stationary distribution}

We introduce an additional constant drift $(-\mu,-\nu)$ in the nonnegative quadrant and consider the system
\begin{equation}
 \label{eq:derive}
 \begin{array}{lll}
 X_t & = & X_0 + B_t + \alpha\int_0^t\frac{ds}{X_s} + \beta \int_0^t\frac{ds}{Y_s} - \mu t\\
 Y_t & = & Y_0 + C_t + \gamma\int_0^t\frac{ds}{X_s} + \delta \int_0^t\frac{ds}{Y_s} - \nu t
 \end{array}
\end{equation}
with the conditions $X_t\geq 0$, $Y_t\geq 0$. If now 
 \begin{equation}
 \begin{array}{lll}
 U_t & = & X_0 + B_t + \alpha\int_0^t\frac{ds}{U_s}   - \mu t\\
 V_t & = & Y_0 + C_t + \delta \int_0^t\frac{ds}{V_s} - \nu t 
 \end{array}
\end{equation} 
with $U_t\geq 0$, $V_t\geq 0$, we can check that ${\bf 0}$ is still polar for $(U,V)$ 
(as well, $0$ is polar for $U$ if $\alpha\geq \frac{1}{2}$ and for $V$ if $\delta\geq \frac{1}{2}$).
Therefore  the results of the previous sections are still valid for the solution to (\ref{eq:derive}). 
We are now looking for conditions on the set 
of parameters in order to obtain a stationary distribution for the Markov process $(X,Y)$ in the form of a product of two
gamma distributions.

\begin{Th}
Assume there exists a unique solution  to (\ref{eq:derive}) in $S^{{\bf 0}}$ or in $S$. This process has 
an invariant distribution in the form $\Gamma(a,c)\otimes \Gamma(b,d)$ if and only if
\[
 \begin{array}{lll}
  \bullet & \alpha\beta+\gamma \delta =0 & \mbox{(skew symmetry)} \\
  \bullet & a=2\alpha+1, & b=2\delta+1 \\
  \bullet & c=2\delta\frac{\mu \alpha +\nu \gamma}{\alpha \delta-\beta\gamma},
     & d=2\alpha\frac{\mu \beta+\nu \delta}{\alpha \delta-\beta\gamma} \\
  \bullet & \mu\alpha+\nu\gamma >0, & \mu \beta+\nu \delta >0 \,.
 \end{array}
\]
\end{Th} 
 {\bf Proof}. Let 
\[
  \rho(x,y)= x^{a-1}e^{-cx} y^{b-1}e^{-dy} \qquad \mbox{for } x\geq 0,y\geq 0 \,.
\]
 The infinitesimal generator of the diffusion (\ref{eq:derive}) is given by
\[
  L= \frac{1}{2}(\frac{\partial^2}{\partial x^2}+\frac{\partial^2}{\partial y^2})
    +(\frac{\alpha}{x}+\frac{\beta}{y}-\mu)\frac{\partial}{\partial x}
    +(\frac{\gamma}{x}+\frac{\delta}{y}-\nu)\frac{\partial}{\partial y} \,.
\]
By a density argument, to prove that $\rho$ is an invariant density,  it is enough to check that 
\[
  \int_0^{\infty}\int_0^{\infty}Lf(x,y) \rho(x,y) \,dxdy\,=0
\]
for any $f(x,y)=g(x)h(y) $ with $g,h \in C^2_c((0,\infty))$ (compactly supported twice continuously 
differentiable functions on $(0,\infty)$). Integrating by parts, we get
\[
  \int_0^{\infty}\int_0^{\infty}L(gh)(x,y)\rho(x,y)\,dxdy\,=\int_0^{\infty}\int_0^{\infty}g(x)h(y)J(x,y)\,dxdy
\]
where 
\[
  J(x,y) = \rho(x,y)\,[A+Bx^{-1}+Cx^{-2}+Dy^{-1}+Ey^{-2}+Fx^{-1}y^{-1}]
\]
with 
\[
  \begin{array}{lll}
  A & = & \frac{1}{2}c^2 + \frac{1}{2}d^2 -\mu c -\nu d \\
  B & = & -(a-1) c + \mu (a-1) + \alpha c + \gamma d \\
  C & = & \frac{1}{2}(a-1)(a-2) -\alpha(a-2) \\
  D & = & -(b-1)d +\nu (b-1) +\beta c + \delta d \\
  E & = & \frac{1}{2}(b-1)(b-2)-\delta (b-2) \\
  F & = & \beta(a-1) + \gamma (b-1)\,.
\end{array}
\]
Letting $A=B=C=D=E=F=0$ we obtain the 
specified values for $a,b,c,d$ and the skew symmetry condition $\alpha\beta+\gamma\delta=0$, which means that 
$r_x$ and $r_y$ are orthogonal. The last condition in the statement of the theorem 
is written out so that the invariant density $\rho$ is integrable on $S$.  It is satisfied if $(\mu,\nu)$ points into the
 interior of the quadrant designed by $r_x$ and $r_y$.  $\hfill \blacksquare$ 
\vspace{0.5cm}

\noindent {\bf Remark.}
With the same proof, we may check that under the skew symmetry condition, 
when $\mu=\nu=0$, the function $\rho(x,y)=x^{2\alpha}y^{2\delta}$ is a non-integrable invariant density that  does not depend 
on the obliqueness parameter $\beta$.
\vspace{0.5cm}

\section{A singular differential system}

In this last section we give up the stochastic setting and consider the integral system
\begin{equation}
 \label{eq:deter}
 \begin{array}{lll}
  x(t)& = & \alpha \int_0^t\frac{ds}{x(s)}\,+\,\beta \int_0^t\frac{ds}{y(s)} \\
  y(t)& = & \gamma \int_0^t\frac{ds}{x(s)}\,+\,\delta \int_0^t\frac{ds}{y(s)}
 \end{array}
\end{equation}
where $x$ and $y$ are continuous functions from $[0,\infty)$ to $[0,\infty)$ with the conditions
\[
\begin{array}{lll}
  \int_0^t{\bf 1}_{\{x(s)=0\}}ds =0 & \quad & \int_0^t{\bf 1}_{\{y(s)=0\}}ds =0\\
  \int_0^t{\bf 1}_{\{x(s)>0\}}\frac{ds}{x(s)}< \infty  & \quad & \int_0^t{\bf 1}_{\{y(s)>0\}}\frac{ds}{y(s)}< \infty 
\end{array}
\]
for any $t\geq 0$.
Here $\alpha,\beta,\gamma,\delta$ are again four real constants with $\alpha>0,\delta>0$ and
the starting point is the corner of the nonnegative quadrant.

An explicit solution is easily obtained.

\begin{Th}
If $\max\{\beta,\gamma\}\geq 0$ or if $\beta\gamma<\alpha\delta$, there is a solution to 
 (\ref{eq:deter}) given by
\[
 \begin{array}{lll}
  x(t) & = & c\sqrt{t} \\
  y(t) & = & d\sqrt{t}
 \end{array}
\]
where
\[
 \begin{array}{lll}
  c & = & (2\alpha+\frac{\beta}{\delta}(\beta-\gamma+\sqrt{(\beta-\gamma)^2+4\alpha\delta}))^{1/2} \\
  d & = & (2\delta+\frac{\gamma}{\alpha}(\gamma-\beta+\sqrt{(\beta-\gamma)^2+4\alpha\delta}))^{1/2} \;.
 \end{array}
\]
There is no solution if $\beta<0,\gamma<0$ and $\alpha\delta\leq \beta\gamma$.
\end{Th}

{\bf Proof}. Writing down $x(t)=c\sqrt{t}$ and $y(t)=d\sqrt{t}$, we have to solve equations
\[
 \begin{array}{lll}
  \frac{c}{2} & = & \frac{\alpha}{c}+\frac{\beta}{d} \\
  \frac{d}{2} & = & \frac{\gamma}{c}+\frac{\delta}{d} 
 \end{array}
\]
Simple computations lead to the requested values. We have to check that
\[
  C=2\alpha+\frac{\beta}{\delta}(\beta-\gamma+\sqrt{(\beta-\gamma)^2+4\alpha\delta})
\]
and
\[
  D=2\delta+\frac{\gamma}{\alpha}(\gamma-\beta+\sqrt{(\beta-\gamma)^2+4\alpha\delta})
\]
are positive. If $\beta\geq 0$, $C$ is clearly positive. This is also true if $\beta<0$ and $\beta\gamma<\alpha\delta$ since $C$ may be written
\[
  C=\frac{4\alpha (\alpha\delta-\beta\gamma)}{2\alpha\delta-\beta\gamma+\beta^2-\beta\sqrt{4(\alpha\delta-\beta\gamma)+(\beta+\gamma)^2}} \;.
\]
The proof for $D$ is similar.

Assume now $\beta<0,\gamma<0,\alpha\delta\leq \beta\gamma$. As was noted in Section 4, there exist $\lambda>0$ and
$\mu>0$ such that $\lambda\alpha+\mu\gamma\leq 0$ and $\mu\delta+\lambda\beta\leq 0$. If $(x,y)$ is a solution to  (\ref{eq:deter}), then   
\[ 
0\leq\lambda x(t)+\mu y(t)\leq 0
\]
  and therefore 
$x(t)=y(t)=0$ for any $t\geq 0$, which is impossible. $\hfill \blacksquare$

We now take up the question of uniqueness.

\begin{Th} Let $\varepsilon= \min
   \{\beta(\beta-\gamma+\sqrt{(\beta-\gamma)^2+4\alpha\delta}),
      \gamma(\gamma-\beta+\sqrt{(\beta-\gamma)^2+4\alpha\delta})\}$.
The solution to  (\ref{eq:deter}) is unique
\[
  \begin{array}{llll}
  a)\mbox{ when } & \beta\geq 0,\gamma\geq 0, & \mbox{if} & \beta\gamma<2\alpha\delta+
    \frac{\varepsilon}{2}\\
  b)\mbox{ when } & \beta\leq 0, \gamma\leq 0, & \mbox{if} & \beta\gamma<\alpha\delta \;.
  \end{array}
\]
\end{Th}

{\bf Proof}. We may use Proposition 2 with $\sigma_j^i=0$ for $i,j=1,2$ and consequently the proofs of Theorem \ref{Th:premier} and 
Theorem \ref{Th:trois} are still valid.

\noindent a) Case $\beta\geq 0, \gamma \geq 0$.\\
As in the proof of Theorem \ref{Th:premier}.3, there exist two extremal solutions $(x_1,y_1)$ and $(x_2,y_2)$ such that for any solution $(x,y)$ and any $t\geq 0$,
\[
  \begin{array}{l}
   x_2(t)\leq x(t)\leq x_1(t) \\
   y_1(t)\leq y(t)\leq y_2(t) \;.
  \end{array}
\]
Uniqueness for equation   (\ref{eq:deter}) follows as well if $\beta\gamma\leq \alpha\delta$. We can improve the result by using 
the explicit form of the known solution and Gronwall lemma. 
Assume $\beta\gamma>\alpha \delta$. As in the proof of Theorem \ref{Th:premier}.2, we consider the unique solution $(x^n,y^n)$ to the system
\[
  \begin{array}{llll}
  x^n(t) & = & \xi_n+\alpha\int_0^t\frac{ds}{x^n(s)} +\beta \int_0^t\frac{ds}{y^n(s)} & \geq 0 \\
  y^n(t) & = & \gamma\int_0^t\frac{ds}{x^n(s)} +\delta \int_0^t\frac{ds}{y^n(s)} &\geq 0  \end{array}
\]
where $(\xi_n)$ is a sequence of decreasing numbers tending to $0$. Remarking that for any $t\geq 0$ we have $c\sqrt{t}\leq  x^n(t)$ and $d\sqrt{t}\geq y^n(t)$, we set
\[
  v_n(t)=(\delta(x^n(t)-c\sqrt{t})-\beta(y^n(t)-d\sqrt{t}))^2 \;.
\]
We develop
\[
  \begin{array}{lll}
    v_n(t) & = & \delta^2\xi_n^2+ 2\int_0^t\sqrt{v_n(s)}(\beta\gamma-\alpha\delta)(\frac{1}{c\sqrt{s}}-\frac{1}{x^n(s)})ds \\
   & = & \delta^2\xi_n^2+2\frac{\beta\gamma-\alpha\delta}{c}\int_0^t\sqrt{v_n(s)}\frac{x^n(s)-c\sqrt{s}}{\sqrt{s}x^n(s)}ds \\
  & \leq & \delta^2\xi_n^2+2\frac{\beta\gamma-\alpha\delta}{c\delta}\int_0^t\frac{v_n(s)}{\sqrt{s}x^n(s)}ds\;.
  \end{array}
\]
From $x^n(t)\geq \xi_n$ and $x^n(t)\geq c\sqrt{t}$ it follows that for $t\geq \xi_n^2/c^2$
$$
  v_n(t) \leq \delta^2\xi_n^2+2\frac{\beta\gamma-\alpha\delta}{c\delta}\left(\int_0^{\xi_n^2/c^2}\frac{v_n(s)}{\xi_n\sqrt{s}}ds+\int_{\xi_n^2/c^2}^t\frac{v_n(s)}{cs}ds\right)\;.
$$
From Gronwall lemma,
\[
  \begin{array}{lll}
   v_n(t) & \leq & \delta^2\xi_n^2\exp\left\{2\frac{\beta\gamma-\alpha\delta}{c\delta}\left(\frac{2\sqrt{\xi_n^2/c^2}}{\xi_n}+\frac{1}{c}\log t -\frac{1}{c} \log \frac{\xi_n^2}{c^2} \right)\right\} \\
  & = & C(t) \xi_n^{2(1-2\frac{\beta\gamma-\alpha\delta}{c^2\delta})}
  \end{array}
\]
where $C(t)$ does not depend on $\xi_n$
  . Therefore $v_n(t)\rightarrow 0$ when $\xi_n\rightarrow 0$ if
\[
  \beta\gamma < 2\alpha\delta+\frac{\beta}{2}(\beta-\gamma+\sqrt{(\beta-\gamma)^2+4\alpha\delta})\;.
\]
We conclude that
\[
  x_1(t)=c\sqrt{t},\;\;y_1(t)=d\sqrt{t}\;.
\]
In a similar way,
\[  
    x_2(t)=c\sqrt{t},\;\;y_2(t)=d\sqrt{t}
\]
under the condition
\[
    \beta\gamma < 2\alpha\delta+\frac{\gamma}{2}(\gamma-\beta+\sqrt{(\beta-\gamma)^2+4\alpha\delta})\;.
\]
Uniqueness follows if both conditions are satisfied.\\
b) Case $\beta\leq 0, \gamma<0$.\\
The proof is the same as for Theorem \ref{Th:trois}.1. $\hfill \blacksquare$
\vspace{0.5cm}

\noindent {\bf Remark.} The question of uniqueness for all $\beta>0,\gamma>0$ remains open. Let us notice that for $\alpha=\delta=0$ uniqueness fails because the solutions are given by 
\[
  \begin{array}{lll}
    x(t) & = & C\,t^{\frac{\beta}{\beta+\gamma}} \\
    y(t) & = & \frac{\beta+\gamma}{C}\,t^{\frac{\gamma}{\beta+\gamma}}
  \end{array}
\]
and depend on the positive parameter $C$.

\vspace{0.5cm}

\footnotesize {{\bf Acknowledgement}. The author is grateful to N. Demni for 
drawing his attention to reference
\cite{OO12}, which motivated this reseach.}

\end{document}